\documentclass[11 pt]{article}
\usepackage[latin1]{inputenc}
\usepackage{amssymb}
\usepackage{hyperref,amsmath,amsfonts,graphicx,psfrag,qsymbols}

\def \pn {\par \noindent}
\def\R{\mathbb{R}}
\def\S{\mathbb{S}}
\def\B{\mathbb{B}}

\DeclareMathOperator{\argmin}{argmin}

\begin{document}
\newtheorem{lem}{Lemma}[section]
\newtheorem{fig}[lem]{\hskip130pt Figure}
\newtheorem{defi}[lem]{Definition}
\newtheorem{theo}[lem]{Theorem}
\newtheorem{fact}[lem]{Fact}
\newtheorem{cor}[lem]{Corollary}
\newtheorem{prop}[lem]{Proposition}
\newtheorem{rem}[lem]{Remark}
\newtheorem{ass}[lem]{Assumption}
\newdimen\AAdi%
\newbox\AAbo%
\font\AAFf=cmex10%   %ou cmex10

\font\dsrom=dsrom10 scaled 1400

\def \cro#1{\llbracket#1\rrbracket}
\def \1{\textrm{\dsrom{1}}}
\def\AArm{\fam0 }%\tenrm}%
\def\eref#1{(\ref{#1})}
\def\AAk#1#2{\setbox\AAbo=\hbox{#2}\AAdi=\wd\AAbo\kern#1\AAdi{}}%
\def\AAr#1#2#3{\setbox\AAbo=\hbox{#2}\AAdi=\ht\AAbo\raise#1\AAdi\hbox{#3}}%
\def \elln#1{\ell_{#1,n}}
\def\ni{\noindent}
\def\proof{\noindent{\bf Proof:}\hskip10pt}
\def\QED{\hfill $\Box$ \smallskip}
\def \be{\begin{eqnarray*}}
\def \ee{\end{eqnarray*}}
\def \ben{\begin{eqnarray}}
\def \een{\end{eqnarray}}

\def \cro#1{\llbracket#1\rrbracket}
\def\d{\displaystyle}
\def\inte{\overset{\circ}}
\def \captionn#1{\begin{center}\begin{minipage}{14cm}\sf\caption{\small #1}\end{minipage}\end{center}}
\def \x{{\sf x}}
\def \w{\xrightarrow[n]{weakly}}
\def \sous#1#2{\mathrel{\mathop{\kern 0pt#1}\limits_{#2}}}
\def \sur#1#2{\mathrel{\mathop{\kern 0pt#1}\limits^{#2}}}
\def\tend {\rightarrow}
\def \dd{\xrightarrow[n]{~(d)~}}
\def \ddk{\xrightarrow[k]{~(d)~}}
\def \as{\xrightarrow[n]{~a.s.~}}
\def \ask{\xrightarrow[k]{~a.s.~}}
\def \pp{\xrightarrow[n]{proba.}}
\def \ppk{\xrightarrow[k]{proba.}}
\def \el{\sur{=}{(d)}}
\def\sn{^{(n)}}
\def\BBa{{\AArm A}}%
\def\BBb{{\AArm I\!B}}%
\def\BBc{{\AArm C\AAk{-1.02}{C}\AAr{.9}{I}{\AAFf\char"3F}}}%
\def\BBd{{\AArm I\!D}}%
\def\BBe{{\AArm I\!E}}%
\def\BBf{{\AArm I\!F}}%
\def\BBg{{\AArm G\AAk{-1.02}{G}\AAr{.9}{I}{\AAFf\char"3F}}}%
\def\BBh{{\AArm I\!H}}%
\def\BBi{{\AArm I\!I}}%
\def\BBj{{\AArm J}}%
\def\BBk{{\AArm I\!K}}%
\def\BBl{{\AArm I\!L}}%
\def\BBm{{\AArm I\!M}}%
\def\BBn{{\AArm I\!N}}%
\def\BBo{{\AArm O\AAk{-1.02}{O}\AAr{.9}{I}{\AAFf\char"3F}}}%
\def\BBp{{\AArm I\!P}}%
\def\BBq{{\AArm Q\AAk{-1.02}{Q}\AAr{.9}{I}{\AAFf\char"3F}}}%
\def\BBr{{\AArm I\!R}}%
\def\BBs{{\AArm S}}%
\def\BBt{{\AArm T\AAk{-.62}{T}T}}%
\def\BBu{{\AArm U\AAk{-1}{U}\AAr{.95}{I}{\AAFf\char"3F}}}%
\def\BBv{{\AArm V}}%
\def\BBw{{\AArm W}}%
\def\BBx{{\AArm X}}%
\def\BBy{{\AArm Y}}%
\def\BBz{{\AArm Z\!\!Z}}%
\def\P{{\AArm I\!P}}%
\def\E{{\AArm I\!E}}%
\def\md{(\Omega ,{\cal I}, {\BBp}^{\theta }, \theta  \in  \Theta )}
\def\BBone{{\AArm 1\AAk{-.8}{I}I}}%
\def\wh{\widehat}
\def\wt{\widetilde}
\def\videbox{\mathbin{\vbox{\hrule\hbox{\vrule height1ex \kern.5em\vrule height1ex}\hrule}}}
\font\calcal=cmsy10 scaled\magstep1
\def\build#1_#2^#3{\mathrel{\mathop{\kern 0pt#1}\limits_{#2}^{#3}}}
\renewcommand{\baselinestretch}{1.2}
\title{On the reduction of a random basis}
\author{Ali Akhavi
\thanks{LIAFA, Université Denis Diderot- Case 7014, 2 place Jussieu, F-75251 Paris Cedex 05 [akhavi@liafa.jussieu.fr]
} \and
Jean-Fran\c cois Marckert
\thanks{LABRI, Université Bordeaux I, 
351 cours de la Libération 
33405-Talence cedex. [marckert@labri.fr] 
}
\and
Alain Rouault
\thanks{LMV, UMR 8100, Universit\'e de Versailles-Saint-Quentin, 45 Avenue des Etats-Unis, 78035-Versailles. [rouault@math.uvsq.fr]}
}

\maketitle
\begin{abstract} 
 For  $g <  n$, let $b_1,\dots,b_{n-g}$ be $n - g$ independent  vectors   in $\mathbb{R}^n$ with a  common distribution invariant by rotation.  Considering these vectors as a basis for the Euclidean lattice they generate, the aim of this paper is to provide asymptotic  results when $n\to +\infty$ concerning the property that such a random basis is reduced in the sense of {\sc Lenstra, Lenstra \& Lov\'asz}. \\
The proof passes by the study of the process $(r_{g+1}^{(n)},r_{g+2}^{(n)},\dots,r_{n-1}^{(n)})$ where $r_j^{(n)}$ is the ratio of lengths of two consecutive vectors $b^*_{n-j+1}$ and $b^*_{n-j}$ built from  $(b_1,\dots,b_{n-g})$ by the Gram--Schmidt orthogonalization procedure, which we believe to be interesting in its own. 
We show that, as $n\to+\infty$, the process $(r_j^{(n)}-1)_j$ tends in distribution in some sense to an explicit process $({\mathcal R}_j -1)_j$; some properties of this latter are provided.
\end{abstract}

\renewcommand{\theequation}{\thesection.\arabic{equation}}
\setcounter{equation}{0}

\section{Introduction.}

We call {\em ambient space} the  space ${\R}^n$ with its classical Euclidean structure. The Euclidean norm is  denoted by $\Vert . \Vert$ and the scalar product by $\langle , \rangle$. Let $b_{1}^{(n)},b_{2}^{(n)},\dots,b_{p}^{(n)}$ (for $p \le n$) be a linearly independent system of $p$  vectors of ${\R}^n$. The superscript $\sn$  is used when needed to stress  the dimension of the ambient space. The quantity 
$$g=n-p,$$
is often used in this paper and referred to as the {\em codimension} of the independent system.

\subsection{The Gram-Schmidt orthogonalization, the reduction level and the index of worst local reduction}\label{GSO}
\pn To the independent system  $b_{1}^{(n)},b_{2}^{(n)},\dots,b_{p}^{(n)}$, the classical Gram-Schmidt orthogonalization procedure associates the orthogonal system 
$\widehat b_1^{(n)} , \cdots , \widehat b_p^{(n)}$ defined by the recursion
\ben
\label{GS}
\widehat b_1^{(n)} = b_1^{(n)} \ , \ \widehat b_j^{(n)} = b_j ^{(n)}- \displaystyle\sum_{i=1}^{j-1}\displaystyle\frac{\langle b_j^{(n)} ,\widehat b_i^{(n)}\rangle}{\Vert \widehat b_i^{(n)}\Vert^2}\ \widehat b_i^{(n)}\,~~~~\textrm{ for }j\geq 2.
\een
If $B = [b_1^{(n)} , \cdots , b_p^{(n)}]$ is the $ n\times p$ matrix with column vectors $b_1^{(n)} , \cdots , b_p^{(n)}$ in the canonical basis, this orthogonalization corresponds to the $QR$ decomposition 
$B = QR$ where 
\[
Q = \left[\widehat b_1^{(n)},\cdots ,  \widehat b_p^{(n)}\right]
\]
is an orthogonal $n\times p$ matrix and $R$ is an upper triangular $p \times p$ matrix $(R_{k, j} = 0 \ , \ 1 \leq j < k \leq n)$ and
\begin{eqnarray}\label{propre}
R_{jj}  =1, \ \ R_{k,j} =  \frac{\langle\widehat b_k^{(n)} , b_j^{(n)}\rangle}{\Vert  \widehat b_k^{(n)}\Vert^2} \ , \ 1 \leq k <j\leq n\, . 
\end{eqnarray}

\begin{defi}\label{defi}

Let $b_{1}^{(n)},b_{2}^{(n)},\dots,b_{p}^{(n)}$  be a linearly independent system  of vectors of $\R^n$ whose  codimension is $g=n-p$. Let $\widehat b_1^{(n)} , \cdots , \widehat b_p^{(n)}$ be the associated Gram-Schmidt orthogonalized system. We call {\em  reduction level} of  $b_{1}^{(n)},b_{2}^{(n)},\dots,b_{p}^{(n)}$  the quantity
\[{\cal M}_n^g := \min_{i\in\{1,\dots,n-(g+1)\}}
 \frac{\Vert\widehat  b_{i+1}^{(n)}\Vert^2}{\Vert\widehat  b_i^{(n)}\Vert^2},
\]
 We call {\em  index of worst local reduction} of $b_{1}^{(n)},b_{2}^{(n)},\dots,b_{p}^{(n)}$ the quantity
\[{\cal I}_n^g := \min \left\{ i : \ \frac{\Vert\widehat  b_{n-i}^{(n)}\Vert^2}{\Vert\widehat  b_{n-i-1}^{(n)}\Vert^2} = {\cal M}_n^g\right\}\,.
\]
\end{defi}

 The motivation of these definitions is explained in Section \ref{S:LLL}. When the vectors $b_{1}^{(n)},b_{2}^{(n)},\dots,b_{p}^{(n)}$ are chosen at random, the reduction level and  the index of worst local reduction are two  random variables, well defined whenever $b_{1}^{(n)},b_{2}^{(n)},\dots,b_{p}^{(n)}$ is a linearly independent system.  This paper is essentially devoted to the study of these random variables. The next subsection details the distribution we consider  for the vectors  $b_{1}^{(n)},b_{2}^{(n)},\dots,b_{p}^{(n)}$.

\subsection{Models of random bases}

In this paper we assume that the $b_i$'s are picked up randomly in $\BBr^n$, independently, and with the same distribution $\nu_n$. Moreover we require $\nu_n$ to be invariant by rotation and to satisfy $\nu_n({0})=0$. 
It is then well known (see \cite{Muir} Th. 1.5.6 p.38 and  Letac \cite{Letaciso}) that the radial part $\Vert b_i^{(n)}\Vert$ and the angular parts 
$\theta_i^{(n)}  := b_i^{(n)}/ \Vert b_i^{(n)}\Vert$ are independent, and that the angular parts are uniformly distributed on ${\S}^{n-1} := \{ x \in \BBr^n : \Vert x\Vert = 1\}$. We call such a model a "simple spherical model".  Since we are interested in the asymptotic behavior of a random basis in $\mathbb{R}^n$ when $n$ goes to $+\infty$, a spherical model will be a sequence of distributions $(\nu_n)$, each $\nu_n$ being a simple spherical model in $\mathbb{R}^n$. \par

The uniform distribution $\mathbb{U}_n$ in the ball ${\B}^n := \{ x \in \BBr^n : \Vert x\Vert \leq 1\}$ -- called the "random ball model" -- is a particular case of spherical model. Under $\mathbb{U}_n$, the distribution of the radial part is  
\ben
\label{rad}
\mathbb{U}_n (\{x : \Vert x\Vert \leq r\}) = \mathbb{U}_n(\Vert b_1^{(n)}\Vert \leq r) = r^n, \ \ 0 \leq r \leq 1\,.
\een 
Under a spherical model, $b_{1}^{(n)},b_{2}^{(n)},\dots,b_{p}^{(n)}$ (for $p \le n$) are a.s. linearly independent. 
We call it a ($p$-dimensional) {\em random basis}.

Our main results hold under assumption (\ref{ass}). This is a technical condition on the distribution $(\nu_n)$ which allows to transfer results concerning the uniform distribution on $\mathbb{S}^{n-1}$ to more general spherical distributions. 

\begin{ass}
\label{ass} There exists a deterministic sequence $(a_n)_n$ and  constants $d_1, d_2 , \alpha>0$, $\rho_0 \in (0,1)$ 
such that, for every $n$ and $\rho \in (0, \rho_0)$
\ben
\label{PGDB}
\nu_n\left(\left|\frac{\Vert b_1^{(n)}\Vert^2}{a_n} - 1\right| \geq \rho\right) \leq d_1 e^{-n d_2 \rho^\alpha}\,.  \een
\end{ass}
This implies in particular that 
$\sup\left\{\left|\frac{\Vert b_i^{(n)}\Vert^2}{a_n} - 1\right|,i\in\{1,\dots,n\}\right\}\xrightarrow[n]{proba}0.$

Here are three natural examples of model $\nu_n$ where such a sequence $(a_n)$ exists: \\
$\bullet$   $\nu_n$ is the uniform distribution on $\mathbb{S}^{n-1}$. In this case $\Vert b_1^{(n)}\Vert^2 = 1$, and $a_n =1$.\\
$\bullet$ $\nu_n=\mathbb{U}_n$. In this case, $a_n=1$ and by \eref{rad}, 
\[\mathbb{U}_n(|{\Vert b_1^{(n)}\Vert^2}/{a_n} - 1| \geq \rho)= (1-\rho)^{n/2}\leq e^{-n\rho/2}.\]\\ 
$\bullet$ $\nu_n$ is the  $n$-variate standard normal (the coordinates are i.i.d. ${\cal N}(0,1)$). Then $\Vert b_1^{(n)}\Vert^2/2$ is $\gamma_{n/2}$-distributed. For $a_n=n$, 
\be
\mathbb{P}(|\Vert b_1^{(n)}\Vert^2/n-1|\geq \rho)&=&\mathbb{P}(\gamma(n/2)\geq (1+\rho)\frac{n}{2})+\mathbb{P}(\gamma(n/2)\leq (1-\rho)\frac{n}{2}).\ee  The Laplace transform $\mathbb{E}(e^{t\gamma(n/2)})$ of $\gamma(n/2)$ is $(1-t)^{-n/2}$, and its Cram\`er transform is 
\ben\label{CT}
H^{(n/2)}(x) = \sup_{\theta < 1}\left\{\theta x-\log \mathbb{E}(e^{\theta\gamma(n/2)})\right\}=x-\frac{n}2+\frac{n}2\log(n/(2x)),~~~ x\geq 0.
\een
 By Markov , $\mathbb{P}(\gamma(n/2)\geq (1+\rho)\frac{n}{2})\leq e^{-H^{(n/2)}((1+\rho)n/2)}=e^{-\frac{n}{2}(\rho-\log(1+\rho))}$ and by an analogous calculus, $\mathbb{P}(\gamma(n/2)\leq (1-\rho)\frac{n}{2})\leq e^{\frac{n}{2}(\rho+\log(1-\rho))}$. Hence assumption \eref{ass} holds in this case with $\alpha=2$. \pn
Notice that these three models are cited in the book of  Knuth (\cite[Section 3.4.1]{Knuth}). \par

 The motivation to study the random variables ${\cal M}_n^g$ and ${\cal I}_n^g$ comes from the theory of ``lattice basis reduction''. The next section briefly describes this motivation and expresses our result in the vocabulary of this theory. The reader who is not interested by this theory may skip the next section.

\subsection{LLL reduction of a random lattice}\label{S:LLL}
 Let $b_{1}^{(n)},b_{2}^{(n)},\dots,b_{p}^{(n)}$ (for $p \le n$) be a linearly independent system of $p$  vectors of ${\R}^n$.
The set of  all their integer linear combinations is an additive discrete subgroup of $\R^n$ called  a lattice. The system $b_{1}^{(n)},b_{2}^{(n)},\dots,b_{p}^{(n)}$ is then a {\em basis} of the lattice. The integer $p$ is the {\em dimension} of the lattice or the dimension of the basis. The codimension  of the lattice basis is the codimension $g=n-p$ of the linearly independent system $b_1^{(n)},\dots,b_p^{(n)}$. The basis  is called {\em full} if  $g=0$.

\medskip

\pn The lattice basis reduction problem  deals with  finding  a basis  of a given lattice,  whose vectors are ``short'' and ``almost orthogonal''.  The problem is old and there are numerous notions of reduction. For a general survey, see for example \cite{Kan87,Val89,Len00}.
Solving even approximately the lattice basis reduction problem  has numerous
theoretical and practical applications in integer optimization \cite{Len83},
computational number theory \cite{LLL82} and cryptography \cite{JoSt98}.

\pn In 1982, Lenstra, Lenstra and Lov\'asz \cite{LLL82}
introduced for the first time   an efficient (polynomial with respect to the length of the input) approximation  reduction algorithm. 
It depends
on a real approximation parameter $s \in ]0, \sqrt{3}/{2}[$ and is called
LLL($s$).  The output basis of the LLL algorithm is called  an LLL($s$) reduced or $s$-reduced basis.  In this paper  we are concerned  with the probability that a random basis under a spherical model is LLL$(s)$ reduced, (i.e. is already an output basis of the LLL($s$)-algorithm). \\
 Roughly speaking the LLL reduction procedure is an approximation algorithm following a divide and conquer paradigm: Indeed for $ i\in\{1\dots p-1\}$, the following condition \eqref{eq:siegel} ensures that some ``local two dimensional basis'' is $s$-reduced. This two dimensional basis is the projections of  $b_i^{(n)}$ and $b_{i+1}^{(n)}$ into the orthogonal $H_i^\perp$ of the vector space $H_i$ spanned by $b_1^{(n)},b_2^{(n)},\dots,b_{i-1}^{(n)}$.  \cite{LLL82} showed  that when all these two--dimensional bases are $s$-reduced then  the whole basis has nice enough Euclidean properties. For instance, the length of the first vector of an LLL-reduced basis is not longer than $(1/s)^{(p-1)}$ times  the length of a shortest vector in the lattice generated by $b_{1}^{(n)},b_{2}^{(n)},\dots,b_{p}^{(n)}$. 

\pn The next definition characterizes an LLL($s$) reduced basis. 

\begin{defi}\label{Siegel}
Let $b_{1}^{(n)},b_{2}^{(n)},\dots,b_{p}^{(n)}$ (for $p \le n$) be a linearly independent system of $p$  vectors of ${\R}^n$. It is an LLL($s$)-reduced basis of the lattice that it generates iff  for all $\ 1 \leq i \leq p-1$, 
\begin{equation}\label{eq:siegel} \frac{\Vert\widehat  b_{i+1}^{(n)}\Vert^2}{\Vert\widehat  b_i^{(n)}\Vert^2} >s^2.
 \end{equation}
\end{defi}
There are two minor differences between the definition of LLL reduction we consider here  and the original definition introduced in \cite {LLL82}.
\pn  Firstly in the original definition the basis has also to be {\em proper}, i.e. if $B= QR$ is the decomposition \eqref{propre} associated with the Gram--Schmidt orthogonalization of the basis  $b_{1}^{(n)},b_{2}^{(n)},\dots,b_{p}^{(n)}$, then 
\begin{equation}\label{e:propre}-1/2 \leq R_{k,j} < 1/2 \ , \ 1 \leq k <j\leq n\, .  \end{equation}

\pn But  from any  basis satisfying \eqref{eq:siegel} one efficiently obtains  a proper basis  still satisfying \eqref{eq:siegel} by a straightforward sequence of integer translations  provided in appendix.  Moreover  considering the notion of {\em flag} \cite{Len00} rather than basis for lattices, makes it possible to skip the notion of properness. 

\pn Secondly the approximation parameter of the original LLL in \cite{LLL82} is slightly different from the one we use here and the reduction we consider here is indeed Siegel reduction as called in \cite{akhavi1,akhavi2}. 
Our main Theorem \ref{T:LLL} is still true with the original definition of a LLL reduced basis as detailed in appendix.

 \medskip

\pn In this paper we study the asymptotics (with respect to the dimension $n$ of the ambient space) of the random variables ${\cal M}_n^g$ and ${\cal I}_n^g$ under spherical models  and for general codimensions of the random basis. 
 The variable ${\cal M}_n^g$ is the supremum of the set of those $s$ for which the basis is $s^2$-reduced. 
As mentioned earlier an LLL($s$) reduced basis  satisfies a set of local conditions. The second variable ${\cal I}_n^g$  is the place where the satisfied local condition is the weakest. This indicates where the limitation of the reduction comes from locally.

\begin{theo}\label{T:LLL}

Let $b_{1}^{(n)},b_{2}^{(n)},\dots,b_{n-g}^{(n)}$ be a random basis with codimension $g$ under  a   spherical model  $(\nu_n)$ satisfying Assumption \eref{ass}. Let $s\in (0,1)$ be a real parameter. 

\pn ($i$)   If  $g=g(n)$ tends to infinity, then the probability
that a random basis is $s$--reduced tends to $1$.
\pn ($ii$) If $g$ is  constant then the probability that a random basis is $s$--reduced converges to a constant in $(0,1)$ (depending on $s$ and $g$). 
\pn ($iii$) If $g$ is  constant, the index of worst local reduction ${\cal I}_n^g$ converges in distribution.
\end{theo}
 Theorem \ref{T:LLL} answers positively to a conjecture of  Akhavi \cite{akhavi1} (which says that for  $c\in[0,1)$, ${\cal M}_n^{cn -1} \pp 1\,$).
In his Lemma 3 p. 376, he proved that
$\P({\cal M}_n^{cn -1} \leq s) \rightarrow 0\,$,
as soon as $s < \frac{1}{2} (1-c)^{\frac{1-c}{c}} (1 +c)^\frac{1}{c}$, and that this convergence is exponentially fast. The proof of Theorem \ref{T:LLL} relies on some properties of random basis under the spherical model which are of interest by their own; these results are overviewed in the next section.

\medskip

\pn Notice that in \cite{Do79}, Donaldson proved a  phenomenon similar to the assertion $(i)$ of Theorem \ref{T:LLL}.  He considered a  different random model:  The basis $b_1^{(n)}, \cdots , b_{n-g}^{(n)}$ is picked up uniformly in the set $\{\Vert b_1^{(n)}\Vert^2 + \cdots + \Vert b_{n-g}^{(n)}\Vert^2 = 1\}$ (Euclidean sphere in $\R^{n\times (n-g)}$). He proved that as $n \rightarrow \infty$ with $n-g(n)$ a fixed constant , the basis is asymptotically reduced in the sense of Minkowski, i.e. each $b_i^{(n)}$ is a shortest vector among all vectors of the lattice that complete $b_1^{(n)}, \cdots , b_{i-1}^{(n)}$ to form a  bigger subset of a lattice basis. So his result is about  a stronger notion of reduction but he considered a much more restricted class of basis. 
\medskip

\pn To finish this Section about lattice basis reduction, observe that our Theorem \ref{T:LLL} about LLL reduction can be generalized to other reductions:  In  \cite{Schn04}  Schnorr introduces a new type of reduction  by segments. In this setting  one fixes an integer $k$ and partitions a basis whose vectors are in $\R^n$ and whose codimension is $g$ into $m$  segments of   $k$ consecutive basis vectors  such that $n-g =km$. For a basis with codimension $g$, the reduction criterion is based on the quantity
\begin{equation}\label{e:schnorr}
M_{k,n}^g = \inf_{r : (k+1)r \leq n-g} \frac{\Vert \wh b^{(n)}_{kr +1}\Vert^2 \cdots \Vert \wh b^{(n)}_{(k+1)r }\Vert^2
}{
\Vert \wh b^{(n)}_{k(r-1) +1}\Vert^2 \cdots \Vert \wh b^{(n)}_{kr }\Vert^2
} \end{equation}
\pn Similarly to the assertions of Theorem \ref{T:LLL}, if $g=g(n)$ tends to infinity and the block size $k$ is fixed, then for any $s\in [0,1]$ the probability that a random basis is $s$--reduced in the sense introduced by Schnorr tends to $1$ with $n$. If $g$ is constant then this probability tends to a constant in $[0,1]$ (depending on $s$, $g$ and $k$).
\footnote{Of course there is a choice of approximation parameters such that  when a basis is LLL($s$) reduced then for any fixed $k$, it is also $s,k$-reduced in the sense introduced by Schnorr. But our approach here shows the existence of limit probabilities (with $n$) for the reduceness of a random basis in the sense introduced by Schnorr.}

\subsection{Random bases issued from spherical models}\label{RL}

For any $j=1,\dots,n$, let 
\[Y_j^{(n)}:=\Vert\widehat  b_j^{(n)}\Vert^2/ \Vert b_j^{(n)}\Vert^2.\]
We denote by $\gamma_a$ and $\beta_{a,b}$ respectively the gamma distribution with parameter $a$, and the beta distribution with parameter $a$ and $b$. In the sequel $\gamma(a)$ and $\beta(a,b)$ stand for generic random variables with respective distribution $\gamma_a$ and $\beta_{a,b}$. Some classical properties of these distributions are recalled in the appendix. \par
We first recall some facts concerning the spherical models, facts that are more or less part of the folklore, and which have been proved several times (e.g.  \cite{Muir},  \cite{akhavi1}).  
\begin{theo}For each $n$, under the simple spherical model, \label{DV}
 the  variables $\Vert \widehat b_j^{(n)} \Vert^2$ ,  $j = 1 , \cdots , n$ are independent.
For every $j=2,\dots,n$,
\ben
\label{Ybeta} Y_j^{(n)} \el \beta\left(\frac{n-j+1}{2}, \frac{j-1}{2}\right)\,,\een
and the random variables $Y_j^{(n)}, j\geq 1$,  $\Vert b_j^{(n)}\Vert^2,j\geq 1$ are independent. 
\end{theo}
A probabilistic proof is given in Section \ref{aerz} for the convenience of the reader.
\begin{cor}
\label{DVorig} Under the random ball model $\mathbb{U}_n$, the variables $\Vert \widehat b_j^{(n)} \Vert^2$, $j = 1 , \cdots , n$ are independent and for $1 \leq j\leq n$ 
\ben
\label{DVorigf}
\Vert\widehat  b_j^{(n)}\Vert^2 \el \beta \left(\frac{n-j+1}{2}, \frac{j+1}{2}\right)\,.\een
\end{cor}
As an easy consequence of the properties of the beta distribution, under $\mathbb{U}_n$,
\ben
\Vert\widehat  b_{n-j}^{(n)}\Vert^2 \el  1 - \Vert\widehat  b_j^{(n)}\Vert^2\,. 
\een 
The statement of Corollary \ref{DVorig} in this formulation is due to Daudé-Vallée (\cite{DD}). Actually, (\ref{DVorigf}) is a consequence of Theorem \ref{DV} and identity (\ref{beta}), since 
(\ref{rad}) means that $\Vert b_i^{(n)}\Vert^2 \el \beta(n/2, 1)$.
 \medskip

\pn The random variable ${\cal M}_n^g$  has the representation :
\ben
\label{rkn}{\cal M}_n^g
=\min_{g+1\leq j\leq n-1} r_j^{(n)}\ ,\ \ r_j^{(n)} := \frac{\Vert\widehat  b_{n-j+1}^{(n)}\Vert^2}{\Vert\widehat  b_{n-j}^{(n)}\Vert^2}
.\een 
As one can guess in view of Theorem \ref{DV}, under $\nu_n$, for each $j$, $r_j^{(n)}$ converges in distribution to $\gamma\left(\frac{j+1}{2}\right)/{\gamma\left(\frac{j}{2}\right)}\,,$
where $\gamma\left(\frac{j+1}{2}\right)$ and $\gamma\left(\frac{j}{2}\right)$ are independent (see Proposition \ref{cvalter}). By the strong law of large numbers, one sees that $\gamma\left(\frac{j+1}{2}\right)/{\gamma\left(\frac{j}{2}\right)}\xrightarrow[j]{a.s.} 1$; this allows to guess that the minimum ${\cal M}_n^g$ is reached by the firsts $r_j^{(n)}$; this motivates the time inversions done in \eref{rkn}.\medskip
\pn The variable ${\cal M}_n^g$ is a function of the $(n-g)$-tuple $(r_{g+1}^{(n)} , \cdots r_{n-1}^{(n)})$, and then the convergence of each coordinate is not sufficient to yield that of ${\cal M}_n^g$. We have to take into account that the variables $(r_j^{(n)})_{j \leq n-1}$ are dependent, and  that their number is growing. Since for the "last" indices ($n-i$ with $i$ fixed),  $r_{n-i}^{(n)} \dd 1$ (see \eref{cvgteralter}), it is convenient to embed the $(n-1)$-tuple $(r_1^{(n)} , \cdots r_{n-1}^{(n)})$ into $\R_+^{\BBn}$ (the set of infinite sequences of positive real numbers), setting 
\ben
\label{embed}
r_j^{(n)} := 1 \ ,\ j \geq n\,.
\een   
Let $(\eta_i)_{i\geq 1}$ be a sequence of independent random variables such that $\eta_i \el \gamma_{i/2}$ and set 
\ben
\label{defr}{\cal R}_j = {\eta_j }/{\eta_{j+1}}\ , \ j \geq 1\,.\een
We denote by $\|.\|_p$ the classical norm on the set  $\ell_p$ of sequences of real numbers with finite $p$th moment: for any sequence of real numbers $\x=(\x_i)_{i \geq 1}$, $\|\x\|_p$ is $\left(\sum_{i\geq 1} |\x_i|^p\right)^{1/p}$ and $\ell_p$ is $\{\x, \|\x\|_p<+\infty\}$.\par
The following result states a limit behavior for the process $(r^{(n)}_j)$ when $n$ goes to$+\infty$.
\begin{prop}\label{lep}For any $p>2$, the following convergence in distribution holds in the metric space $\ell_p$~:
\[(r_j^{n}-1)_{j\geq 1} \xrightarrow[n]{(d)} ({\cal R}_j-1)_{j\geq 1}.\]
\end{prop}
The previous proposition is the key result here and it will entail all the convergence results given in the next theorem.

\pn For $k \in \BBn$, set
\[{\cal M}^k:=\inf\big\{{\cal R}_j, j\geq k+1\big\}.\]
The application ${\sf x} \mapsto 1+\min_{i\geq k} {\sf x}_{i}$ is continuous from $\ell_p$ onto $\mathbb{R}$.
It follows that ${\cal M}_n^g\wedge 1$  converges in distribution to ${\cal M}^g$. We will prove that 
\begin{theo}
\label{t2} If $\nu_n$ is spherical and satisfies Assumption \ref{ass} then,\\ 
$(i)$ For each $k$, ${\cal M}_n^k \dd  {\cal M}^k$.\\
$(ii)$ Let $g:{\mathbb{N}}\tend \mathbb{N}$ such that $g(n)\leq n$ and $ g(n) \rightarrow \infty$. We have ${\cal M}_n^{g(n)}\pp 1\,. $\\
$(iii)$~For any $k\geq 1$, ${\cal I}_n^k\dd{\cal I}^k$.
\end{theo}

\pn Notice that Proposition \ref{lep} and Theorem \ref{t2} have their analogous for the reduction introduced by Schnorr in \cite{Schn04}. By setting
$$
\label{rtidekn}{\cal M}_{k,n}^g
=\min_{g+1\leq kr\leq n-1} r_{k,r}^{(n)}\ ,\ \ r_{k,r}^{(n)} := \frac{\Vert\widehat  b_{n-(r+1)k+1}^{(n)}\Vert^2 \dots \Vert\widehat  b_{n-rk}^{(n)}\Vert^2}{\Vert\widehat  b_{n-(r+2)k+1}^{(n)}\Vert^2 \dots \Vert\widehat  b_{n-(r-1)k}^{(n)}\Vert^2}
\text{ and }\ \ r_{k,r}^{(n)} := 1 \text{ for }kr\geq n,
$$ 
if we let $n \rightarrow \infty$, we have convergence of  $(r_{k,r}^{(n)})_r$ to a process $({\cal R}_{k,r})_r$ with
$${\cal R}_{k,r} = \frac{\eta_{k, r}}{\eta_{k,r+1}} \ , \ \eta_{k,r} \el \gamma(r/2) \gamma((r+1)/2\cdots \gamma((r+k-1)/2)\,,$$ 
where the $\eta_{k,r}, r \geq 1$ are independent, and the gamma variables too. Then by setting 
\[\widetilde{{\cal M}}_k^g:=\inf\big\{{\cal R}_{k,r}, kr\geq g+1\big\}.\]
 one obtains also an analogous to  Theorem \ref{t2}.

\medskip

\begin{figure}[htbp]
\centerline{
\includegraphics[height=5cm,width =7.5cm]{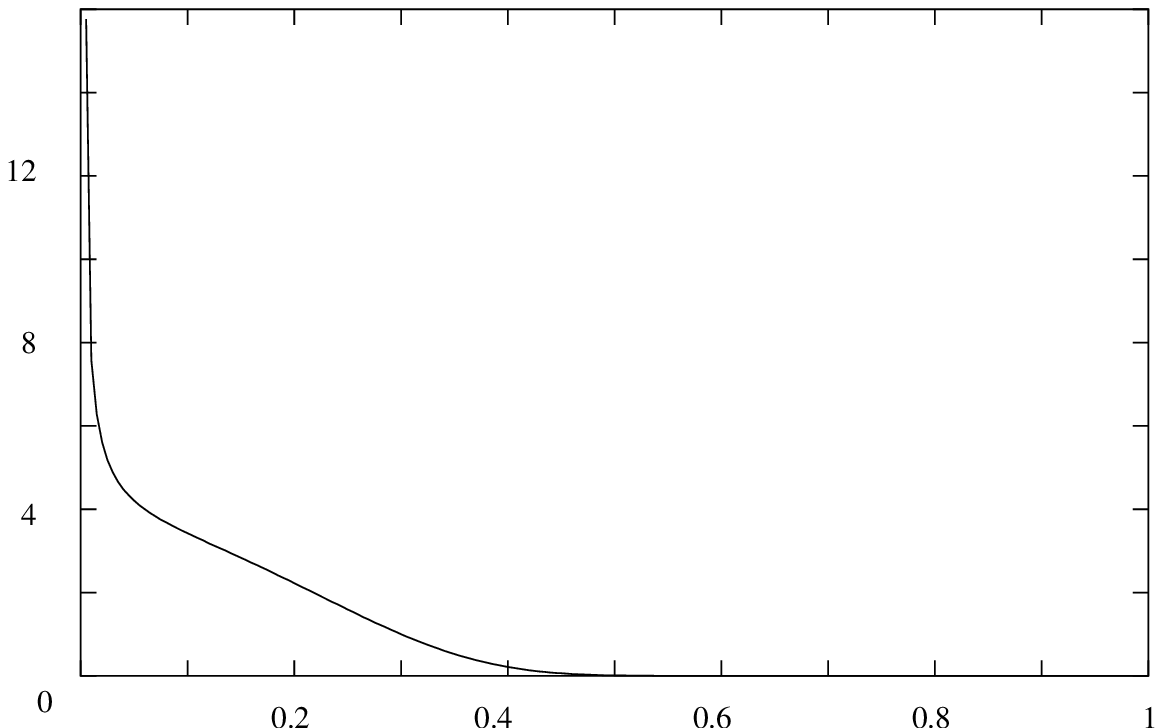}\includegraphics[height=5cm, width =7.5 cm]{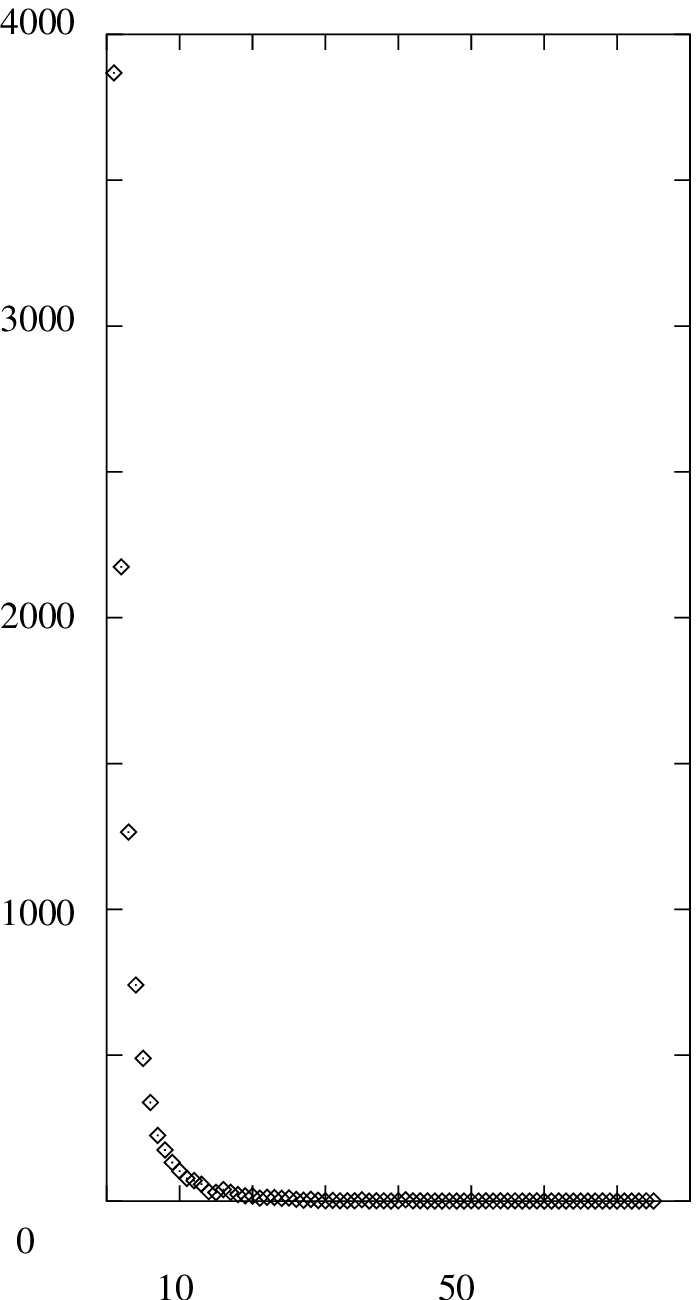}}
\captionn{\label{sim2} On the first picture, simulation of the density of ${\cal M}_{\infty}^0$ with $10^8$ data. On the second, the histogram provided by 10000 simulations of ${\cal I}_{\infty}$. The sequence $k\mapsto P({\cal I}^g=k)$ seems to be decreasing.} 
\end{figure}

\noindent \bf A note on the proof of Proposition \ref{lep}. \rm  The ambient spaces $\R^n, n \geq 1$ are not nested, and then we give up the geometrical consideration on $\mathbb{R}^n$ and focus on the representation of the processes $r_j^n$ using the gamma distributions.  \medskip

We end this section by stating some properties of the limiting process $({\cal R}_k)_{k\geq 1}$. First of all, in statistics the distribution of $\frac{j+1}{j}{\cal R}_j$  is known as the Fisher $F_{j, j+1}$-distribution (its distribution is recall in (\ref{quotient})); the mean of ${\cal R}_j$ is $j/(j-1)$ and, as said above, ${\cal R}_k \ask 1$. Here are some sharper results (see also simulations on Figure \ref{sim2}).

\begin{prop}\label{density} 
$(i)$
For  each $k$, the distribution of ${\cal M}^k$ has a density, which is positive on $(0,1)$ and zero outside.\\
 \label{elen}
$(ii)$ For each $k$,
\[ \lim_{x \downarrow 0} x^{-(k+1)/2} \P({\cal M}^k\leq x) = {1}/{\Gamma\left(\frac{k+2}{2}\right)}.\]
$(iii)$ There exists $\tau>0$ such that for each $k\geq 0$, 
\[\limsup_{y \uparrow 1} \!\ e^{\frac{\tau}{(1-y)^2}}\!\ \P({\cal M}^k\geq y) < \infty\,.\]
$(iv)$ For each $k$,  there is a.s. a unique random index ${\cal I}^k$ such that ${\cal R}_{{\cal I}^k} = {\cal M}^k$.
\end{prop}

\section{Proofs} 
\label{proofs}
\setcounter{equation}{0}
\subsection{Additional information on random basis}
\label{aerz}

We first give a proof of Theorem \ref{DV} for convenience.\medskip

\noindent{ \bf Proof of Theorem \ref{DV}}: Let us skip the superscript $(n)$ in this proof. We have $b_i = \theta_i \Vert b_i\Vert$ and from (\ref{GS}), we see that $\wh b_i = \Vert b_i\Vert \wh\theta_i$, where the $\wh\theta_i$'s are obtained by the Gram-Schmidt algorithm applied to the $\theta_i$'s. The independence of $(\wh\theta_1, \cdots, \wh\theta_n)$ and $(\Vert b_1\Vert^2, \cdots , \Vert b_n\Vert^2)$ is then a direct consequence of the radial-angular independence. Notice that $\wh\theta_1 = 1$. 
 
Now, fix $k\geq 2$. Conditionally upon $\wh\theta_1, \cdots, \wh\theta_{k-1}$, the variable $\Vert \wh\theta_{k}\Vert$ is  distributed as the norm of the projection of a random vector uniformly distributed on ${\mathbb{S}}^{n-1}$  on $span \{\wh\theta_1, \cdots, \wh\theta_{k-1}\}$. 
Since the problem is invariant by rotation, this distribution is independent of $(\wh\theta_1, \cdots, \wh\theta_{k-1})$ which proves (recursively) that  the $\Vert\wh\theta_i\Vert$'s are independent. 
 Moreover,   $\Vert\wh\theta_k\Vert$ is distributed as the norm of the projection of $\theta_k$ (or $\theta_1$) on the subspace generated by the $n-k+1$ 
last vectors of the canonical basis. From  Muirhead \cite[Theorem 1.5.7, p. 38-39]{Muir} the distribution of $\Vert\wh\theta_k\Vert^2$ is $\beta_{\frac{n-k+1}{2}, \frac{k-1}{2}}$.
\QED

Here are some information on the asymptotic behavior of the random variables $Y_j^{(n)}$~:
\begin{prop}Under a spherical model, \label{alter}
for each $j\geq 1$,
\ben
\label{cvalter}
\frac{n}{2} \   Y_{n-j}^{(n)} &\dd& \gamma_{\frac{j+ 1}{2}}\,,\\
\label{cvgbisalter}
 Y_{j}^{(n)}&\dd& 1\,.
\een
\end{prop}

In view of Theorem \ref{DV} this yields:
\begin{prop}Under a spherical model, \label{t1alter}
if $\Vert b_1^{(n)}\Vert^2/a_n \dd 1$ for some deterministic sequence $a_n$, then
for each $j\geq 1$,
\ben
\label{cvgalter}
\frac{n}{2a_n} \  \Vert\widehat  b_{n-j}^{(n)}\Vert^2 &\dd& \gamma_{\frac{j+ 1}{2}}\,,\\
\label{cvgteralter}
\frac{1}{a_n}\ \Vert\widehat  b_{j}^{(n)}\Vert^2 &\dd& 1\,.
\een
\end{prop}

\begin{rem}
Under the same assumptions, we have also:

If $h(n) \rightarrow \infty$ and   $h(n)/n \rightarrow 0$,  then
\begin{equation}
\label{cvg1}
\frac{n}{h(n)a_n} \ \Vert\widehat  b_{n-h(n)}^{(n)}\Vert^2 \pp 1 \ .
\end{equation}
If $0<\alpha < 1$ et $k(n)/n \rightarrow 0$, then
\begin{equation}
\label{cvg2}\frac{1}{a_n}\ \Vert\widehat  b_{\alpha n + k(n)}^{(n)}\Vert^2 \pp 1 - \alpha \  .
\end{equation}
\end{rem}

This result stated under $\mathbb{U}_n$ can be found in \cite[Theorem 8]{akhavi1}. Let us give a new proof which prefigures the main 
arguments used to prove the convergences in Section \ref{PC}.
\medskip

\noindent{\bf Proof of Propositions \ref{alter} and \ref{t1alter}}
From Theorem \ref{DV} we have the decomposition,
\ben
\label{y}
\Vert\widehat  b_{n-j}^{(n)}\Vert^2 \sur{=}{(d)} Y_{n-j}^{(n)}\ \Vert b_1^{(n)}\Vert^2\,,
\een with $Y_{n-j}^{(n)} \sur{=}{(d)} \beta\left(\frac{j+1}{2}, \frac{n-j-1}{2}\right)\,.$
Let $(\xi_j )_{j\geq 1}$ be a sequence of i.i.d. $\gamma_{1/2}$-distributed  random variables. From (\ref{gamma1}) and (\ref{gamma0}), we can write
\ben
\label{eqrep}
 Y_{n-j}^{(n)} \sur{=}{(d)} \displaystyle\frac{\sum_{m= 1} ^{j+1} \xi_m}{\sum_{m= 1} ^n \xi_m}.
\een
By the strong law of large numbers, 
\[\frac{\sum_{m= 1} ^n \xi_m}{n} \xrightarrow[n]{a.s.} \frac{1}{2},\] and for each $j$, $\sum_{m= 1} ^{j+1} \xi_m \sur{=}{(d)} \gamma((j+1)/2)$, which yields \eref{cvalter}. 
From this and the additional assumption $\Vert b_1^{(n)}\Vert^2/a_n \dd 1$, we see that (\ref{cvgalter}) holds true. For \eref{cvgbisalter}, notice that 
$(1 - Y_j^{(n)}) \sur{=}{(d)} Y_{n-j+2}^{(n)}\,,$ and that $Y_{n-j+2}^{(n)}\pp 0$ by \eref{cvalter}. To end,  (\ref{cvgteralter}) is a consequence of \eref{cvgbisalter} and $\Vert b_1^{(n)}\Vert^2/a_n \sur{\to}{(d)} 1$. \QED 
\medskip

The following lemma will be used to transfer results from the uniform distribution on $\mathbb{S}^{n-1}$ to more general spherical distributions. 
\begin{lem}
\label{isocalc}
Assume that Assumption \ref{ass} holds.  If $U_1$ and $U_2$ be independent and $U_1 \el U_2 \el \Vert b_1^{(n)}\Vert^2$, then there exist $d'_1, d'2, \alpha > 0$ and $\rho_0 \in $ such that for any $k \geq 1, n \geq 1$ and $\rho\in (0, \rho_0)$ 
\ben
\P \left(\left|\frac{U_1}{U_2} - 1\right|\geq \rho\right) \leq d'_1 \exp(- nd'_2 \rho^\alpha\,).
\een
\end{lem}

\proof 
We have
\be
\P \left(\frac{U_1}{U_2} \geq 1 + \rho\right) &\leq& \P (U_2 \leq (1-\rho/2)) + \P(U_1 \geq (1+\rho)(1-\rho/2))\\ 
&\leq&\P (U_2 \leq (1-\rho/2)) + \P(U_1 \geq (1+\rho/4)))
\ee
as soon as $\rho \leq 1/2$. 
Similarly
\be
\P \left(\frac{U_1}{U_2} \leq 1 - \rho\right) &\leq& \P (U_2 \geq (1+\rho/2)) + \P(U_1 \leq (1-\rho)(1+\rho/2))\\ 
&\leq&\P (U_2 \geq (1+\rho/2)) + \P(U_1 \leq (1-\rho/2))
\ee
With the help of assumption (\ref{PGDB}), this yields 
\be
\P \left(\left|\frac{U_1}{U_2} - 1\right|\geq \rho\right) \leq d'_1 \exp(-n d'_2 \rho^\alpha\,).
\ee
\QED

\subsection{The process $({\cal R}_k)$: estimates and proof of Proposition \ref{density}}
Lemma \ref{CLT} and Proposition \ref{LDP} first state some properties concerning the fluctuations and large deviations of the distribution ${\cal R}_k$ 
\begin{lem} The following convergence in distribution holds
\label{CLT}
$$\sqrt{k}\!\ ({\cal R}_k - 1)\ddk {\cal N}(0, 4)\,.$$
\end{lem}
\bf Proof : \rm 
Setting \[\xi_k = \frac{\eta_k - k/2}{\sqrt k} ~~\textrm{ and }~~ \xi'_k = \frac{\eta_{k+1} - (k+1)/2}{\sqrt k}\]
the CLT gives $(\xi_k , \xi'_k)\xrightarrow[k]{~~(d)~~}  {\cal N} (0, 1/2)\otimes {\cal N} (0, 1/2)$ hence $\xi_k  - \xi'_k \xrightarrow[k]{~~(d)~~} {\cal N}(0, 1)$.
Since
\[\sqrt{k}\!\ ({\cal R}_k - 1) = \frac{k}{\eta_{k+1}} \left(\xi_k  - \xi'_k -\frac{1}{2\sqrt k}\right) 
\,, \]
and $\eta_{k+1}/k \rightarrow 1/2$ a.s.,
 we get the result. \QED 
\begin{prop}
\label{LDP}
Let $f_{{\cal R}_k}$ be the density of ${\cal R}_k$ and \be
\Phi_k (x) =
{(4x)^{\frac{k}{2}-1}}(1 +x)^{-k - \frac{1}{2}}\,.
\ee 
\begin{enumerate}
\item For $A < {2}{ \pi}^{-1/2} < B$ we can find an integer $K$ such that
\ben
\label{dens1}
A \sqrt{k}\ \Phi_k(x)\leq f_{{\cal R}_k} (x)\leq B \sqrt{k} \ \Phi_k (x)
\een
for every $x\in (0, \infty)$ and every $k \geq K$.
\item There exists a constant $C$ such that for every $k \geq 1$ and $\rho \in [0,1]$
\ben
\label{majo2}
\P({\cal R}_k < 1 - \rho) 
  &\leq& C \left( 1 - \frac{\rho^2}{(2 - \rho)^2}\right)^{k/2}\\
\label{majo1}
\P({\cal R}_k > 1 + \rho) 
 &\leq& C \left( 1 - \frac{\rho^2}{(2 + \rho)^2}\right)^{k/2}\,.
\een
\item Assertion 2 holds true when ${\cal R}_k$ is replaced by ${\cal R}'_k := \displaystyle\frac{S_k^{(2)}}{S_k^{(1)}}$\,.
\end{enumerate}
\end{prop}
Notice that the distribution of ${\cal R}'_k$ is known in statistics as the Fisher $F_{k,k}$.\par
\proof 
1) We have $f_{{\cal R}_k} (x) = C_k  \Phi_k (x)$ where
\be
C_k = 4^{1 -\frac{k}{2}}
 \frac{\Gamma\big(k + \frac{1}{2}\big)}{\Gamma\big(\frac{k}{2}\big)\Gamma\big(\frac{k+1}{2}\big)}
= \frac{2}{\sqrt{\pi }}\frac {\Gamma\big(k + \frac{1}{2}\big)}{\Gamma(k)}
\sim \sqrt{k} \frac{2}{\sqrt \pi}
\,.
\ee
2) The bounds may be obtained by integration, but 
 also by
 writing the beta variables as ratios of gamma variables and using Chernov's bounds.
Noticing that ${\cal R}_k$ and ${\cal R}'_k$ are Fisher-distributed, the above results are related
to section 4 of \cite{Cha}.
Since we need bounds holding for $\rho$ depending on $k$, we use the classical Chernov's method :
\be
\P({\cal R}_k  > 1 + \rho) &=& \P \left(\eta_k - (1 + \rho)\eta_{k+1} > 0 \right)\\
&\leq& E \exp \left( \theta\eta_k 
- \theta (1 + \rho) \eta_{k+1}  \right)
= \left(E e^{\theta \eta_1 }\right)^k \left(E e^{- \theta (1 + \rho)\eta_1 }\right)^{k+1}\\
&=& (1- \theta )^{-k/2} \left(1 + \theta (1 +\rho)\right)^{-(k+1)/2}\\
&=& \left(1 + \theta (1 +\rho)\right)^{-1/2} \left((1-\theta) (1 + \theta(1 +\rho)\right)^{-k/2}\,. 
\ee
The function $\theta\mapsto (1-\theta) (1 + \theta(1 +\rho))$ reaches its maximum for
$\theta = \frac{\rho}{2(1+ \rho)}\in (0, 1)$, so that :
\begin{equation}
\P({\cal R}_k  > 1 + \rho ) \leq
\left( 1 - \frac{\rho^2}{(2 + \rho)^2}\right)^{k/2}\,.
\end{equation}
Similarly
\be
\P ({\cal R}_k < 1 - \rho ) &\leq& E\exp \left( \theta(1-\rho)
\eta_{k+1} - \theta 
\eta_k \right)\\
&=& \left((1+\theta))(1 - \theta(1 -\rho)\right)^{-k/2} \left(1 - \theta (1 -\rho)\right)^{-1/2}\\
&\leq& \sqrt{2} \left( 1 - \frac{\rho^2}{(2 + \rho)^2}\right)^{k/2}\,.
\ee
3) For ${\cal R}'_k$ the proof needs similar evaluations and is left to the reader.
\QED

Thanks to these bounds on the deviation of the process $({\cal R}_k)$ around the value 1, one may establish the following corollary.
\begin{cor}\label{lp}
For any $p > 2$, the process $({\cal R}_k-1)$ is a.s. in $\ell_p$, i.e.
$\sum_k |{\cal R}_k - 1|^p < \infty$ a.s..
\end{cor}
\bf Proof \rm 
Thanks to the Borel-Cantelli lemma,
it is enough to find
 $v=(v_k)_{k\geq 1} \in \ell_p$,
such that
\ben
\label{aprouver}
\sum_k \P (|{\cal R}_k - 1| \geq v_k) < \infty\,.
\een
Taking $\rho = k^{-\mu}$ in the bounds (\ref{majo1}) and (\ref{majo2}),  we have
$\sum_k \P\left( |{\cal R}_k  - 1 | > k^{-\mu}\right) < \infty$
if $1-2\mu > 0$.
For $p>2$, one may choose $\mu\in]1/p,1/2[$ and $v_k=k^{-\mu}$. Then $(v_k)_{k\geq 1} \in \ell_p$ and satisfies (\ref{aprouver}). \QED

\subsubsection*{Proof of Proposition \ref{density}}
\bf Proof of $(i)$. \rm  We give a proof in the case $k=0$, but the  argument is 
the same  for any $k>0$.\par
First, since  for any $j$, ${\cal R}_j> 0$ a.s. and since a.s., $\lim_j {\cal R}_j= 1$, the support of ${\cal M}^0$ is included in $[0,1]$. For the same reason, the sequence $({\cal R}_k)$ does not accumulate at $0$, which yields that the distribution of ${\cal M}^0$ has no atom at $0$. \par
Using Lemma \ref{CLT} write 
\[\P({\cal R}_j  < 1) = \P\left(\sqrt{2j}\!\ ({\cal R}_j - 1) < 0\right)\xrightarrow[j\rightarrow \infty]{}{1}/{2}.\]
Hence by the reverse Borel Cantelli lemma, a.s. there exists an infinite sequence of $j$ such that ${\cal R}_{2j} < 1$, which yields that ${\cal M}^0$ has no atom at $1$.\par
It remains to check that the support of ${\cal M}^0$ is exactly $[0,1]$ (see (1) below) and that ${\cal M}^0$ has a density (see (2) below).\\
(1)  
Let us prove that $\P(\inf_j {\cal R}_j \in [a,b]) > 0$, for every $[a,b] \subset [0,1]$. It is enough to
find a sequence of (independent) events  $B_j := \{\eta_j \in (\alpha_j , \beta_j)\}, j \geq 0$
such that
\ben
\label{claim}\bigcap_{j=1}^\infty B_j \subset \left\{\inf_j {\cal R}_j \in [a,b]\right\} ~~~\textrm{ and }  ~~~
 \prod_{j=1}^\infty \P(B_j) > 0.
 \een
Let $j_a = \inf\{j : j > 2(1+a)/(1-a)\}$, $A := j_a (1+a)/4$ and $c_1< c_2$ in $(a,b)$. Choose 
\ben
\nonumber
\alpha_1 =Ac_1,\ \alpha_2 = A&,& \ \alpha_j = A \ \hbox{for}\ 3 \leq j \leq j_a\ \ , 
\ \alpha_j = \frac{j(1+a)}{4}  \ \ \  \ \ \  \ \ \ \hbox{for}\ j \geq j_a +1\,,\\
\nonumber
\beta_1 = Ac_2,\ \beta_2 = \frac{Ac_1}{a}&,& \ \beta_j = \frac{A}{a} \ \hbox{for}\ 3 \leq j \leq j_a\ \ , \ \beta_j = \frac{(j-1)(1+a)}{4a}  \ \hbox{for}\ j \geq j_a +1\,.
\een

We check easily that  $B_1 \cap B_2 \subset \{{\cal R}_ 1 \in (a, c_2)\}$,
and $B_j \cap B_{j+1} \subset \{{\cal R}_j \in (a , \infty)\}$
for $j \geq 2$. This proves the first claim of  (\ref{claim}).

 It remains to prove that the infinite product is convergent, i.e. that
\ben
\label{prodinfini}\sum_{k > j_a} \P(B_k^c) < \infty\,.
\een
For $j > j_a$, the interval $(\alpha_j , \beta_j)$ straddles the mean $j/2$ of $\eta_j$ :
$$\alpha_j = \frac{j (1+a)}{4} < \frac{j}{2} \ \ ,\ \  \beta_{j} \geq \frac{j(1+3a)}{8a} > \frac{j}{2}\,,$$
so that the large deviations inequalities hold:
\[\log \P(\eta_j < \alpha_j) \leq
 -j H^{(1/2)}\left(\frac{1+a}{4}\right)\ , \ \log \P(\eta_j > \beta_{j}) \leq  -j H^{(1/2)}\left(\frac{1+3a}{8a}\right)\]
where $H^{(1/2)}$, the Cram\'er transform of $\gamma_{1/2}$ is given in \eref{CT}. This yields a positive constant $M$  such that for $j> j_a$
\[\P(B_j^c) = \P(\eta_j < \alpha_j) + \P(\eta_j > \beta_{j}) \leq 2 e^{-jM}\]
and the series is
convergent, which proves (\ref{prodinfini}) and $\P(\inf_j {\cal R}_j \in [a,b]) > 0$.
\medskip

(2)  According to Radon-Nikodym's theorem,  
it suffices to find a positive integrable function $f$ on
$(0,1)$, such that for any $[a,b]\subset(0,1)$,
\[\P({\cal M}^0\in [a,b])\leq \int_{[a,b]}f(x)dx\]
By the union bound, we have for every $b' \in (b, 1)$ :
\[ \P({\cal M}^0 \in[a,b]) = \P(\inf_{k\geq 1} {\cal R}_k\in[a,b])\leq \P\left(\cup_k \{{\cal R}_k \in [a , b'
]\}\right)
\leq \sum_{k\geq 1}\P\big(
{\cal R}_k \in[a,b']\big)\,.\]
 For $B>2/\sqrt{\pi}$, thanks to formula (\ref{dens1}), there exists $K \geq 1$ such that 
\be
\sum_{k\geq K} \P\big({\cal R}_k \in[a,b]\big) &\leq& B \int_a ^{b'} \left(\sum_{k\geq K} \sqrt{k}\ \Phi_k(x)\right)\!\ dx
\\
&\leq& \frac{B}{2} \int_a ^{b'} \left[\sum_{k\geq 1
} k \left(\frac{2\sqrt x}{1+x}\right)^{k-1}\right]\ \frac{dx}{\sqrt x (1+x)^{3/2}} \\
&=& \frac{B}{2} \int_a ^{b'} \frac{\sqrt{1+x}}{\sqrt x (1- \sqrt x)^4}\ dx\,.
\ee
Since every ${\cal R}_k$ has a density, one may bound the $K-1$ first terms of the sum by $\int_a^{b'} f_1(x) dx$ for 
some integrable $f_1$. Then, since the bound holds true for any $b'>b$, 
we can let $b'\downarrow b$ and we get the result. 

\bf Proof of $(ii)$ \rm  We have $\P({\cal R}_{k+1}\leq x)\leq \P({\cal M}^k\leq x)\leq \sum_{j\geq k+1}\P({\cal R}_j\leq x).$
Using (\ref{quotient}), one obtains, for $x \rightarrow 0$,
\[\P({\cal R}_{k+1}\leq x)=\frac{\Gamma\left(\frac{k+3}{2}\right)\,\,x^{(k+1)/2}}{\frac{(k+1)}2\Gamma\left(\frac{k+1}2\right)\Gamma\left(\frac{k+2}2\right)}\ (1 + o(1)) = \frac{x^{(k+1)/2}}{\Gamma\left(\frac{k+2}2\right)} \ (1 + o(1))\,.\]
On the other hand, a simple computation shows that, when $x \rightarrow 0$,
\[\sum_{j\geq k+2}\P( {\cal R}_j\leq x)=O(x^{(k+3)/2}).\]
\bf Proof of $(iii)$ \rm  We have, for $j \geq k$
\be
\P\left({\cal M}^k> 1 - j^{-1/2}
\right) \leq \prod_{i=j}^{2j} \P\left({\cal R}_{2i} > 1 -
j^{-1/2}
\right)
\leq \prod_{i=j}^{2j} \P \left({\cal R}_{2i} > 1 - i^{-1/2}
\right)\,.
\ee
From Lemma \ref{CLT}, we know that
$\lim_k \P\left({\cal R}_{2k} > 1 -k^{-1/2}
\right) = \P (N > - \sqrt{2})$
where $N$ is ${\cal N} (0, 4)$.
Taking $\tau >0$ with $e^{-\tau} > \P(N > - \sqrt{2})$ we see that for $j$ large enough
\be
\P\left({\cal M}^k > 1 - j^{-1/2}
\right) \leq e^{-\tau j}
\ee
which ends the proof of $(iii)$. \\
\bf Proof of $(iv)$ \rm The support of ${\cal M}^k$ is $[0,1]$ and $\lim {\cal R}_j = 1$  a.s. so that the set $\{j\geq k+1,  {\cal R}_j={\cal M}^k\}$ is not empty. Moreover since there are no ties ($\P({\cal R}_i = {\cal R}_j) = 0$ a.s. for $i\not=j$) this set is a.s. a singleton.\QED

\subsection{The proofs of convergence (Theorem \ref{t2} and Proposition \ref{lep})}
\label{PC}

In order to prove Proposition \ref{lep} and Theorem \ref{t2}, we build a probability space on which are defined some copies of the variables $\|b^{(n)}_i\|$, $i\geq 0$, $n\geq 0$ (and then also $r_j^{(n)}$) and the process $({\cal R}_k)$. This space is not related with some embedding of $\mathbb{R}^n$ in some larger space:  the proof is not geometrical.  Thanks to that procedure, we will be able to use the strong law of large numbers obtaining in such a way strong versions of the convergences in distribution stated in Proposition \ref{lep} and Theorem \ref{t2}.\medskip 

From Theorem \ref{DV} and the representation (\ref{gamma1})
 we see that 
\ben\nonumber
 \Vert \wh b^{(n)}_{n-k+1}\Vert^2  &=& Y_{n-k+1}\sn \Vert b_{n-k+1}\sn\Vert^2\\ Y_{n-k+1}\sn \el \displaystyle\frac{\sum_{m= 1} ^{k} \xi_m}{\sum_{m= 1} ^{n} \xi_m}\ &,& \  \Vert b_{n-k+1}\sn\Vert^2 \el \Vert b_{1}\sn\Vert^2\een
where the $\xi_m$'s are $\gamma_{1/2}$ distributed, and  $\Vert b_{n-k+1}\sn\Vert^2$ is independent of the $\xi_m$'s.
 Since the  $ \Vert \wh b^{(n)}_{n-k+1}\Vert^2$ for  $1\leq k\leq n-1$ are independent, we may consider two double arrays $(\xi_i^k, i \geq 1, k\geq 1)$, $(\zeta_j^k, j \geq 1, k\geq 1)$ of independent random variables (and independent together), such that

a) for every $j \geq 1$ and $k \geq 1$, $\xi_j^k \el\gamma(1/2)$, 

b) for every $j \geq 1$ and $k \geq 1$, $\zeta_j^k \el \Vert b_1^{(j)}\Vert^2$.\par
The common probability space on which are defined all the variables $\xi_{j}^k$ and $\zeta_j^k$ is denoted by $\Omega$. From now on we work exclusively on $\Omega$. \medskip  

Let us set
$$S_{p}^k = \sum_{m=1}^p \xi_{m}^k,\ \ k\geq 1, \ \ p\geq 1\,.$$
Now, the processes $(S_j^{k})_{j\geq 1}$ for $k=1, \cdots$ are independent copies of $(S_j^1)_{j\geq 1}$, and for each $n\geq 1$,
we have the following distributional representation :
\begin{equation}\label{qq}
\{ 
\Vert \wh b^{(n)}_{n-k+1}\Vert^2\   ,  \ \ 1\leq k\leq n-1 \}
\sur{=}{(d)} \left\{ \frac{S_k ^{k}}{S_{n} ^{k}}\ \zeta_n^k  ,  \ \ 1\leq k\leq n-1 \right\}\,.
\end{equation}
For $n\geq 2$, set
\ben
\label{defrkn}
 R_k ^{(n)}&=&
     \begin{cases}
       \frac{\displaystyle S_{k} ^{k} S_{n} ^{k+1}}{\displaystyle S_{k+1}^{k+1} S_{n}^{k}}\  \displaystyle\frac{\zeta_n^k}{\zeta_n^{k+1}} & \text{if $1 \leq k \leq n-1$} \\
         1 & \text{if $k \geq n$}
     \end{cases}
\een
we have now,  (see \eref{rkn} and \eref{embed})
\ben
\label{el} r^{(n)}
\sur{=}{(d)}R^{(n)}\,.\een

The processes $r^{(n)}, n \geq 2$ are not defined on a unique probability space, 
since the ambient spaces are not nested.
On the contrary, the sequence $R^{(n)} , n\geq 2$ is  defined on the unique probability space $\Omega$.
For each $k\geq 1$, the strong law of large numbers yields
\[ \frac{S_{n}^{k+1}}{n} \as \frac{1}{2}
\ , \ \frac{S_{n}^{k}}{n} \as \frac{1}{2}\,,
\]
Besides, Lemma \ref{isocalc}, with the help of Borel-Cantelli's lemma yields 
 \[\frac{\zeta_{n}^{k}}{\zeta_{n}^{k+1}} \as 1\,,\]
so that if we set
\ben \label{eaz}
{\cal R}_k:= \frac{ S_{k} ^{k}}{S_{k+1}^{k+1}}\,,
\een
we get for any $k\geq 1$ 
\[R_k ^n \as  {\cal R}_k.\] Notice that  ${\cal R}$  is defined in \eref{defr}. Hence, letting ${\cal R}_k$ be $\frac{ S_{k} ^{k}}{S_{k+1}^{k+1}}$ here is a slight abuse of notation but this is consistent in terms of distribution and allows to avoid a new symbols. From now on $({\cal R}_k)$ is then a random variable on $\Omega$.
 Setting, for any $g\geq 0$,
\begin{equation}
M_n^g=\min_{g+1\leq k\leq n-1} R_k^{(n)}  ~~~\textrm{ and }~~~{\cal M}^g=\min_{k\geq g+1} {\cal R}_k\,,
\end{equation}
we get
\begin{equation}\label{arko}
{\cal M}_n^g~~~\sur{=}{(d)}~~~M_n^g 
\,,\end{equation}
and want to prove a convergence (in probability) of $M_n^g$ to ${\cal M}^g$. 
Since the convergence of the coordinates of $R^{(n)}$ to those of $({\cal R}_k)$ is not sufficient to this aim, we need
 a uniform control.

Set
\[\widetilde {M}_n^g := \inf_{k\geq g+1}  R_k\sn ,\]
so that $\widetilde {M}_n^g = {M}_n^g \wedge 1$. This yields  $0\leq {M}_n^g - \widetilde {M}_n^g = ({M}_n^g -1)^+ \leq (R_{n-1} ^n - 1)^+$.
Since $ R_{n-1}\sn \pp 1$ (by Theorem \ref{t1alter}), we get
\begin{equation}\label{iin}
{M}_n^g - \widetilde {M}_n^g \pp  0,
\end{equation}
and so, ${M}_n^g$ and $\widetilde {M}_n^g$ have the same limit behavior.\par

To prove Theorem \ref{t2}, we first assume that the following lemma which is a strong form of Proposition \ref{lep} holds true
\begin{lem}
\label{l1}
For any $p > 2$, $( R_k\sn -{\cal R}_k)$ converge a.s. (in $\Omega$) to 0 in $\ell_p$, i.e. 
\begin{equation}\label{cpabossa}
\sum_{k=1}^\infty| R_k\sn -{\cal R}_k|^p \as 0.
\end{equation}
\end{lem}
\smallskip

\subsubsection*{Proof of Theorem \ref{t2}}

 $(i)$ \rm From (\ref{cpabossa}) and Lemma \ref{lp}, the sequence $\big(R_k\sn - 1\big)_{ k\geq 1}$ converges
a.s.
 in $\ell_p$ to
$\left({\cal R}_k  - 1\right)_{k\geq 1}$.  Let $K$ be a fixed integer.
Since the mapping $(c_k)_{k\geq 1} \in \ell_p \longmapsto \inf_{k \geq K} c_k$
is continuous, one has
\begin{equation}\label{resuu}
\widetilde {M}_n^K
~~~\as~~ {\cal M}^K.
\end{equation}
Thanks to \eref{iin}, we obtain
${M}_n^K  \pp  {\cal M}^K$ and then, by \eref{arko} ${\cal M}_n^K\dd {\cal M}^K$.

$(ii)$  Let $`e>0$ and $`e'>0$ be fixed.
Since  $({\cal R}_k-1)_{k\geq 1}\in \ell_p$,  there exists $K$ such that
\[\P({\cal M}^K\leq 1- `e/2)\leq `e'.\]
For $n$ large enough, one then has, by \eref{resuu} and \eref{iin},
\[\P({M}_{n}^K\leq 1-`e)\leq 2`e'.\]Since the function $k\mapsto M_n^k$ is non-decreasing, one has, for $n$ large enough such that $g(n)\geq K$,
\[\P({M}_{n}^{g(n)}\leq 1-`e)\leq 2`e'.\]
$(iii)$ Take $k= 0$ for the sake of simplicity. For $a \in \R^{\BBn}$, let $\argmin a = \{i : \inf_{j \geq 1} a_j = a_i\}$ and as usual set $\min \emptyset = \infty$.
Denote by $I_n^0=\min\argmin\{ R^n_j,j\geq 1\}$,
the natural version of ${\cal I}_n^0$  on $\Omega$~:
\begin{equation}\label{marcelle}
I_n^0\sur{=}{(d)}{\cal I}_n^0
\end{equation}We know that a.s. ${\cal M}^0 < 1$ so that for $n$ large enough, we have 
${M}_n^0 < 1$, hence 
 $\argmin \widetilde R^n = \argmin R^n$.
Now, from Proposition \ref{l1}$(ii)$, a.s $\lim \widetilde R^n = {\cal R}$ in $\ell_p$. 
Now, the convergence of $y_n$ to $y$ in $\ell_p$ implies the convergence of
$\min\argmin(y_n)$ to $\argmin(y)$ if $\#\argmin(y)=1$.
Hence, a.s. $\lim I_n^0 = {\cal I}^0$. Thanks to \eref{marcelle},  
we deduce ${\cal I}_n^0\dd {\cal I}^0$.\QED

\subsubsection*{Proof of Lemma \ref{l1}.} 
Set $V_n := \sum_k |R_k\sn - {\cal R}_k|^p =
  V'_n + V''_n $ where
$$V' _n := \sum_{1 \leq k\leq n-1} |R_k\sn - {\cal R}_k|^p \ \hbox{and} \
V'' _n := \sum_{k \geq n} |1 - {\cal R}_k|^p .$$
According to Lemma \ref{lp}, $V''_n \as 0$.
 Then, it is enough to prove that  $V' _n \as 0$.
Since
\[R_k\sn =  {\cal R}_k \frac{S_{n}^{k+1}}{S_{n}^{k}}\ \frac{\zeta_n^{k}}{\zeta_n^{k+1}}\]
 and since
$\sup_{k \geq 1}{\cal R}_k$ is a.s. finite, it is enough to prove that
\ben
\label{true}
\sum_{k=1}^{n-1} \left|\frac{S_{n} ^{k+1}}{S_{n}^{k}}\ \frac{\zeta_n^{k}}{\zeta_n^{k+1}} - 1 \right|^p
 \as 0\,.
\een
Let $\delta > 0$. By the union bound and the identity of distributions, we have
\ben
\nonumber
\P\left(\sum_{k=1}^{n-1} \left|\frac{S_{n} ^{k+1}}{S_{n}^{k}}\ \frac{\zeta_n^{k}}{\zeta_n^{k+1}} - 1 \right|^p > 
\delta \right) &\leq& \sum_{k=1}^{n-1}
\P\left( \left|\frac{S_{n} ^{k+1}}{S_{n}^{k}}\ \frac{\zeta_n^{k}}{\zeta_n^{k+1}} - 1 \right|^p > \frac{\delta}{n} \right)\\
\nonumber
&=& (n-1) \P\left( \left|\frac{S_{n}^{(2)}}{S_{n}^{(1)}}\ \frac{\zeta_n^{1}}{\zeta_n^{2}} - 1 \right| > \frac{\delta^{1/p}}{n^{1/p}} \right)\,.
\een
Splitting this event, we get easily for $\varepsilon = \frac{\delta^{1/p}}{n^{1/p}}$
\[
\P\left( \left|\frac{S_{n}^{(2)}}{S_{n}^{(1)}}\ \frac{\zeta_n^{k}}{\zeta_n^{k+1}} - 1 \right| > \varepsilon \right)\leq \P\left( \left|\frac{S_{n}^{(2)}}{S_{n}^{(1)}} -1\right| > \varepsilon/3 \right) + \P\left( \left|\frac{\zeta_n^{1}}{\zeta_n^{2}} - 1 \right| > \varepsilon/2\right)
\] 
With the notation of the preliminaries, the first probability  is
$\P(|{\cal R}'_{n} - 1| > \varepsilon/3)$.
By a simple calculation using (\ref{majo1}),(\ref{majo2}) and lemma \ref{isocalc}, we can find  $c_1$ and $c_2 > 0$ such that for every $n$
\be
\P\left(\sum_{k=1}^n \left|\frac{S_{n} ^{k+1}}{S_{n}^{k}}\ \frac{\zeta_n^{k}}{\zeta_n^{k+1}} - 1 \right|^p > \delta \right) 
\leq c_1 n \exp\left(-c_2 n^{1-\frac{2}{p}}\right)\,.
\ee
For $p > 2$, we get  a convergent series, so (\ref{true}) holds true, which ends the proof of $ii)$.\QED

\section{Appendix}
\subsection{LLL($\delta$)-reduced basis versus Siegel($s$)-reduced basis}
As mentioned in Section \ref{S:LLL}, the definition \ref{Siegel} is slightly different from the original definition of an LLL reduced basis as defined in \cite{LLL82}. Here we make precise this difference and show that our main result (Theorem \ref{T:LLL})  is still true with the original definition.
 
Let $(b) := b_{1}^{(n)},b_{2}^{(n)},\dots,b_{p}^{(n)}$ (for $p \le n$) be a linearly independent system of $p$  vectors of ${\R}^n$ and recall the definition of the matrix $R$ given in Section \ref{GSO}.
\begin{defi}\label{LLL} 
Let $0<\delta<1$ be a real parameter.  The basis $(b)$ is called truly--LLL($\delta$) reduced if it is proper \eqref{e:propre} and if 
\begin{equation}\label{eq:LLL}\forall i\in \{ 1\dots,n-1\},\quad 
{\Vert\widehat  b_{i+1}\Vert^2} + R^2_{i,i+1}{\Vert\widehat  b_i\Vert^2} > \delta^2 {\Vert\widehat  b_i\Vert^2}.
\end{equation}
\end{defi}
\pn From the above definition and the definition of a LLL($s$)-reduced basis \eqref{Siegel}, and since $R^2_{i+1,i} \leq 1/4$ (thanks to the properness) one deduces immediately:
\begin{fact}
\begin{itemize}
\item[($i$)] If a basis is LLL($s$) reduced and proper then it is truly--LLL($s$) reduced.
\item[($ii$)] If a basis is truly--LLL($\delta$) reduced then it is LLL($\sqrt{\delta^2 -1/4 }$) reduced.
\end{itemize}
\end{fact}

\subsection{How to make a basis proper while preserving its LLL reduceness}
\pn Here is a simple enunciation of the  LLL$(\delta)$ algorithm:
\medskip
\pn {\bf The Make--proper  algorithm}:
\vspace{1mm}
\pn {\bf Input:} A basis $b=(\bold{b}_1,\dots,\bold{b}_p)$ of a lattice $L$.
\pn {\bf Output:} A  proper basis $b$ of the lattice 
$L$. 

\pn {\bf Initialization:} Compute the orthogonalized system $\widehat{b}$ 
and the matrix $R$.
\pn {\bf For i from 2 to n do }
\par {\bf For j from (i-1) downto 1 do }
\par \hspace{1cm}$\bold{b}_{i}:=\bold{b}_{i}-\lfloor R_{j,i} \rceil \bold{b}_j$ ($\lfloor x 
\rceil$
is the integer nearest to $x$).
\medskip
\pn Clearly the Gram-Schmidt basis associated with the input basis is preserved under  the integer translations of the above algorithm. So the Gram Schmidt orthogonalized basis associated with the output basis is the same as the one associated with the input basis and the   Make--proper  algorithm preserves LLL($s$)-reduceness and truly--LLL($s$)-reduceness. 

\medskip

\subsection{A brief description of the LLL algorithm}

\pn In this subsection, we provide  a simple enunciation of the  LLL$(\delta)$ algorithm. Clearly if the input basis is LLL($s$)-reduced and proper then it is also truly LLL($s$)-reduced. So in this case the following algorithm will stop after one iteration of the \bf while loop \rm (which makes the basis proper).

\medskip
\pn {\bf The $\mathbf{LLL(\delta)}$-reduction algorithm}:
\vspace{1mm}
\pn {\bf Input:} A basis $b=(\bold{b}_1,\dots,\bold{b}_n)$ of a lattice $L$.
\pn {\bf Output:} A $LLL(\delta)$-reduced basis $b$ (or a truly LLL($s$)-reduced basis) of the lattice 
$L$. 

\pn {\bf Initialization:} Compute the orthogonalized system $\widehat{b}$ 
and the matrix $R$.
\pn $\mathbf{i:=1;}$
\pn {\bf While $\mathbf{i<n}$ do}
\par $\bold{b}_{i+1}:=\bold{b}_{i+1}-\lfloor R_{i,i+1} \rceil \bold{b}_i$ ($\lfloor x 
\rceil$
is the integer nearest to $x$).
\par {\bf Test: } ${\Vert\widehat  b_{i+1}\Vert} > s {\Vert\widehat  b_i\Vert}$ ? (or ${\Vert\widehat  b_{i+1}\Vert^2} + R^2_{i,i+1}{\Vert\widehat  b_i\Vert^2} > \delta {\Vert\widehat  b_i\Vert^2}$ ?)

\par \hspace{5mm} {\bf If true,} make $(\bold{b}_{1},\dots, \bold{b}_{i+1})$ 
proper by {\bf Make-proper}; 
{\bf set }$\mathbf{ i:=i+1;}$
\par \hspace{5mm} {\bf If false,} swap $\bold{b}_{i}$ and $\bold{b}_{i+1}$; update 
$\widehat{b}$ and $R$; if $i\neq 1$ then {\bf set }$\mathbf{ i:=i-1;}$
\medskip

\subsection{The Beta--Gamma algebra}
\setcounter{equation}{0}
We recall some properties of the Gamma and Beta distribution, used all along the lines of the paper. They can be found in \cite{ChauYor} pp. 93-94.
For $a>0$, the gamma distribution of parameter $a$ is
\[\gamma_a(dx) = \frac{e^{-x} x^{a- 1}}{\Gamma (a)} \ \BBone_{[0, \infty)} (x)\ dx\,,\]
and its mean is $a$.

For $(a,b)\in {\mathbb R}^{+\star}$, the beta distribution of  parameters $(a, b)$ denoted  by $\beta_{a,b}$ is
$$
\beta_{a,b} (dx) = \frac{\Gamma (a+b)}{\Gamma (a) \Gamma(b)} x^{a-1} (1-x)^{b-1}
\ \BBone_{(0,1)} (x) \ dx.
$$

In the following, $\gamma(a)$ denotes a  variable with distribution $\gamma_a$, and $\beta(a,b)$ denotes a  variable with distribution $\beta_{a,b}$. The first relation is
\ben
\big(\gamma(a) , \gamma(b) \big)&\sur{=}{(d)}&\big(\beta(a,b)\gamma(a+b) , (1 - \beta(a,b)) \gamma(a+b)\big)\,,
\een
where, on the  left hand side the random variables $\gamma(a)$ and $\gamma(b)$ are independent and on the right  hand side the random variables $\beta(a,b)$ and $\gamma(a+b)$ are independent.
It entails
\ben
\label{gamma0}
\gamma(a) + \gamma(b) \sur{=}{(d)}\gamma (a+b) \,,
\een
\ben
\label{gamma1}
\frac{\gamma(a)}{\gamma(a) + \gamma(b)}&\sur{=}{(d)}& \beta(a,b)\,,
\een
and
\ben
\label{gamma2}
\frac{\gamma (a)}{\gamma (b)} &\sur{=}{(d)}& \frac{\beta(a,b)}{1 - \beta(a,b)}\,,
\een
which gives
\ben
\label{quotient}
\P\big({\gamma (a)}/{\gamma (b)} \in dx\big) =
\frac{\Gamma(a+b)}{\Gamma(a)\Gamma(b)}\frac{x^{a-1}}{(1+x)^{a+b}} \ \BBone_{[0, \infty[} (x)\ dx\,.
\een
The second relation is
\ben
\label{beta}
\beta(a,b)\beta(c, a-c) &\sur{=}{(d)}& \beta(c, a+b-c)\,,
\een
where on the left hand side the random variables are independent.

\small\renewcommand{\baselinestretch}{1}

\bibliographystyle{plain}
\bibliography{alinew4}

\begin{thebibliography}{10}

\bibitem{akhavi2}
A.~Akhavi.
\newblock {\em Analyse comparative d'algorithmes de r\'eduction sur les
  r\'eseaux al\'eatoires}.
\newblock PhD thesis, Universit\'e de Caen, 1999.

\bibitem{akhavi1}
A.~Akhavi.
\newblock Random lattices, threshold phenomena and efficient reduction
  algorithms.
\newblock {\em Theoretical Computer Science}, 287:359--385, 2002.

\bibitem{Cha}
N.R. Chaganthy.
\newblock Large deviations for joint distributions and statistical
  applications.
\newblock {\em Sankhya}, 59:147--166, 1997.

\bibitem{ChauYor}
L.~Chaumont and M.~Yor.
\newblock {\em Exercises in probability}.
\newblock Cambridge Series in Statistical and Probabilistic Mathematics.
  Cambridge University Press, Cambridge, 2003.
\newblock A guided tour from measure theory to random processes, via
  conditioning.

\bibitem{DD}
H.~Daudé and B.~Vallée.
\newblock {An upper bound on the average number of iterations of the LLL
  algorithm.}
\newblock {\em Theor. Comput. Sci.}, 123(1):95--115, 1994.

\bibitem{Do79}
J.L. Donaldson.
\newblock Minkowski reduction of integral matrices.
\newblock {\em Mathematics of Computation}, 33(145):201--216, 1979.

\bibitem{Len00}
Jr. H.W.~Lenstra.
\newblock Flags and lattice basis reduction.
\newblock In {\em European Congress of Mathematics, Vol. I (Barcelona, 2000)},
  volume 201 of {\em Progr. Math.}, pages 37--51. Birkh\"auser, Basel, 2001.

\bibitem{Kan87}
R.~Kannan.
\newblock Algorithmic geometry of numbers.
\newblock In {\em Annual review of computer science, Vol.\ 2}, pages 231--267.
  Annual Reviews, Palo Alto, CA, 1987.

\bibitem{Knuth}
D.~E. Knuth.
\newblock {\em The art of computer programming. {V}ol. 2}.
\newblock Addison-Wesley Publishing Co., Reading, Mass., second edition, 1981.
\newblock Seminumerical algorithms, Addison-Wesley Series in Computer Science
  and Information Processing.

\bibitem{LLL82}
A.~K. Lenstra, H.~W. Lenstra, Jr., and L.~Lov{\'a}sz.
\newblock Factoring polynomials with rational coefficients.
\newblock {\em Math. Ann.}, 261(4):515--534, 1982.

\bibitem{Len83}
H.~W. Lenstra, Jr.
\newblock Integer programming and cryptography.
\newblock {\em Math. Intelligencer}, 6(3):14--19, 1984.

\bibitem{Letaciso}
G.~Letac.
\newblock Isotropy and sphericity: some characterisations of the normal
  distribution.
\newblock {\em Ann. Statist.}, 9(2):408--417, 1981.

\bibitem{Muir}
R.~J. Muirhead.
\newblock {\em Aspects of multivariate statistical theory}.
\newblock John Wiley, 1982.

\bibitem{JoSt98}
P.~Q. Nguyen and J.~Stern.
\newblock The two faces of lattices in cryptology.
\newblock In {\em Cryptography and lattices (Providence, RI, 2001)}, volume
  2146 of {\em Lecture Notes in Comput. Sci.}, pages 146--180. Springer, 2001.

\bibitem{Schn04}
C.P. Schnorr.
\newblock {Fast {LLL}-Type Lattice Reduction}.
\newblock {\em Information and Computation}, 204:1--25, 2006.

\bibitem{Val89}
B.~Vall{\'e}e.
\newblock Un probl\`eme central en g\'eometrie algorithmique des nombres: la
  r\'eduction des r\'eseaux. {A}utour de l'algorithme de {L}enstra {L}enstra
  {L}ovasz.
\newblock In {\em Informatique Théorique et Applications}, volume~3, pages
  345--376. 1989.
\newblock English translation by E. Kranakis CWI-Quaterly - 1990 - 3.

\end{thebibliography}

\end{document}